# ANALYTICAL RAMIFICATIONS OF DERIVATIVES VALUATION: ASIAN OPTIONS AND SPECIAL FUNCTIONS

Michael Schröder

*Lehrstuhl Mathematik III*
*Seminargebäude A5, Universität Mannheim, D–68131 Mannheim*

Averaging problems are ubiquitous in Finance with the valuation of the so–called Asian options on arithmetic averages as their most conspicuous form. There is an abundance of numerical work on them, and their stochastic structure has been extensively studied by Yor and his school. However, the analytical structure of these problems is largely unstudied. Our philosophy now is that such valuation problems should be considered as an extension of the theory of special functions: they lead to new problems about new classes of special functions which should be studied in terms of and using of the methods of special functions and their theory. This is exemplified by deriving integral representations for the Black–Scholes prices based on Yor's Laplace transform ansatz to their valuation. They are obtained by analytic Laplace inversion using complex analytic methods. The analysis ultimately rests on the gamma function which in this sense is found to be at the base of Asian options. The results improve on those of Yor and have served us a as starting point for deriving first time benchmark prices for these options.

**1. Introduction:** This paper studies analytical ramifcations of derivatives valuation which it considers as an extension of the theory of special functions. More precisely, the focus is on the analytical aspects of averaging problems. These are ubiquitous in Finance and make their appearance in as diverse areas as interest rate derivatives, credit risk, and problems originating in life insurance contracts. We consider their most conspicuous form: the so–called Asian options on arithmetic–averages. These are widely traded options whose prices, in a way made precise in §4, are quoted using mathematical models. There is by now an abundance of numerical work on how to compute these prices with almost daily additions like [**V**]. However, even in the rather simplifying but still basic Black–Scholes model reviewed in §3, the conceptual understanding of their valuation appears to be far from simple and has lead early researchers to conclude that 'it is impossible to derive an explicit analytic expression for an Asian option price'. And as it can be seen from the discussion in [**RS**] for instance, coming to grips with valuing Asian options has intrigued financial theorists for over a decade by now.

An explanation for this has in a sense been found by Yor and his collaborators whose work has defined the state of the art on this problem in particular. Establishing connections with Bessel and confluent hypergeometric functions they have studied the stochastic structure of Asian option valuation with [**Y01**] as the lastest addition. With regrad to closed form valuation formulas for Asian option prices, Yor defined in [**Y**] the standard with his triple integral. And as a second line of thought, the Laplace transform approach developed in [**GY**] was a celebrated advance. Its results, and its actual relation to Asian option valuation, is reviewed in §8 and furnish the starting point of our research.

The main mathematical result of this paper is a second valuation integral for the Black–Scholes price of Asian options which improves on Yor's. As described in §6 it expresses this price as a single integral over now, roughly speaking, the product of two well–studied higher transcendental functions. Both are given as integrals and built up using so–called Hermite functions, reviewed in §2. These last functions generalize the familiar complementary error function and come from boundary–value problems in potential theory for domains whose surface is an infinite parabolic cylinder. It seems to be the first time that such parabolic cylinder functions are explicitly identified to characterize solutions to problems of Finance. Moreover, our results improve on those difficulties with numbers of gigantic sizes which render Yor's triple integral so uncomputable. This we discuss in §7. And, in fact, our formula improves on these difficulties to such an extent that it has served in [**SH**] and [**ST**] as the starting point for deriving first time explicit benchmark pricing formulas for Asian options.

Mathematically, our valution formula is obtained by analytic Laplace inversion of certain of Yor's Laplace transforms. And this can be regarded as the main mathematical contribution of the paper. The details of the proof are described in Part III of this paper. Its methods are complex analytic. It uses the complex inversion formula for the Laplace transform and is ultimately based on the classical Hankel formulas for the gamma function. However, inversion is achieved under restrictions only that are of stochastics origin. These are then lifted in a second step using analytic continuaton. This seems to be a new feature which illustrates how complex analytic methods can be used to extend the validity of stochastics results in situations where stochastic methods alone seem to lead to intractable problems.

And all this illustrates part of our guiding philosophy that such valuation problems seem best considered as extensions of the theory of special functions. They give rise to new problems in and new classes of special functions which should be studied in terms of and using the methods of special functions and their theory.

**Acknowledgements:** It is a pleasure to thank my friends J. Ballmann, D. Fulea, and U. Weselmann, now of Universität Heidelberg, for their support and the genuine interest they have taken in this whole project, and I would like to include hereby also Dr. P. Carr (Courant Institute, New York) and Professors Yor (Paris VI) and Pliska (UIC, Chicago). Finally, I am very grateful to Professor Harder for making possible my very pleasant stay at the *Mathematisches Institut der Universität Bonn* in the *Wintersemester 97/98* during which work for this paper was started.

# Part I Preliminaries and recollections

**2. Preliminaries on Hermite functions:** This section collects relevant facts about Hermite functions from [**L**, §§10.2ff]. These functions are given in terms of the Kummer confluent hypergeometric function $\Phi$, and the *Hermite function* $H_\mu$ of degree any complex number $\mu$ is the function on the complex plane defined by

$$H_\mu(z) = \frac{2^\mu \cdot \Gamma(1/2)}{\Gamma\big((1-\mu)/2\big)} \cdot \Phi\Big(-\frac{\mu}{2}, \frac{1}{2}, z^2\Big) + z \cdot \frac{2^\mu \cdot \Gamma(-1/2)}{\Gamma\big(-\mu/2\big)} \cdot \Phi\Big(\frac{1-\mu}{2}, \frac{3}{2}, z^2\Big).$$

Hermite functions $H_\mu$ so specialize to the $\mu$–th Hermite polynomials if $\mu$ is any non–negative integer. Moreover, they are holomorphic on the complex plane as functions of both their variable $z$ and their degree $\mu$. If $\mu$ is not a non–negative integer, they have the absolutely and compactly convergent series

$$H_\mu(z) = \frac{1}{2 \cdot \Gamma(-\mu)} \sum_{n=0}^{\infty} \frac{(-1)^n}{n!} \, \Gamma\left(\frac{n-\mu}{2}\right) (2z)^n,$$

which also describes their behaviour for small arguments. If the real part $\mathrm{Re}\,(\mu)$ of $\mu$ is negative, the asymptotics of $H_\mu$ for large arguments is described by the expansion

$$H_\mu(z) = (2z)^\mu \sum_{k=0}^{n-1} \frac{(-\mu)_{2k}}{k!} \cdot \frac{(-1)^k}{(2z)^{2k}} + O\Big(\frac{1}{|z|^{2n-\mathrm{Re}\,(\mu)}}\Big)$$

valid for any complex $z$ with positive real part. Herein recall $(\alpha)_k$ as the Pochhammer symbol given by $(\alpha)_0 = 1$ and $(\alpha)_{k+1} = (\alpha+k)(\alpha)_k$ for $k$ any non–negative integer.

To discuss further connections of Hermite functions with other classes of special functions, if the real part of $\mu$ is negative, we have the integral representation:

$$H_\mu(z) = \frac{1}{\Gamma(-\mu)} \int_0^\infty e^{-u^2-2zu} u^{-(\mu+1)}\, du\,,$$

for any complex number $z$. Hermite functions thus also generalize the *complementary error function* Erfc which for any complex number $z$ is given by:

$$\mathrm{Erfc}\,(z) = e^{-z^2} \frac{2}{\sqrt{\pi}} \int_0^\infty e^{-u^2-2zu} du\,;$$

in fact, $(2/\sqrt{\pi})H_{-1}(z) = \exp(z^2)\mathrm{Erfc}\,(z)$. Hermite functions of degree any complex number $\mu$ are moreover connected with the parabolic cylinder functions $D_\mu$ by $D_\mu(z) = 2^{-\mu/2}\exp(-z^2/4)H_\mu(2^{1/2}z)$ for any complex $z$, and with the second Kummer confluent hypergeometric function $\Psi$ by $H_\mu(z) = 2^\mu\Psi(-\mu/2, 1/2; z^2)$ if $\mathrm{Re}\,(z) > 0$.

**3. Black–Scholes modelling:** The analytical problems to be discussed originate from the so–called risk–neutral approach to the valuation of contingent claims. Textbook treatments for this and other notions developed for the analysis of financial markets and instruments are in [**DSM**], [**MR**], [**KSb**, Chapters 1–4] for instance. This analysis is based on modellings of security markets, and this section aims to give some background on the most fundamental of these, the Black–Scholes modelling of security markets.

In fact, we need only that particular case of the Black–Scholes model where there are only two securities., and the understanding is is that these are traded on markets where their prices are determined by equating demand and supply. First, there is a riskless security, a bond, whose price $\beta$ grows at the continuously compounding positive interest rate $r$, i.e., for which we have we have $\beta_t = \exp(rt)$ at any time $t \in [0,\infty)$. Then, there is a risky security, and the fundamental idea is that all uncertainties affecting its price $S$ yield a certain probability space. In fact, consider for this a complete probability space equipped with the standard filtration of a standard Brownian motion on the time set $[0,\infty)$. Giving expression to the fact that $S$ comes as an equilibrium price, we have the risk neutral measure $Q$ on this filtered space, a probability measure equivalent to the given one. And with $B$ any standard $Q$–Brownian motion, the exact modelling then is that $S$ is the strong solution of the following stochastic differential equation:

$$dS_t = \varpi \cdot S_t \cdot dt + \sigma \cdot S_t \cdot dB_t, \qquad t \in [0,\infty)\,,$$

or equivalently using Itô calculus,

$$S_t = S_0 e^{(\varpi - \frac{1}{2}\sigma^2)t + \sigma B_t} \qquad t \in [0, \infty)\,.$$

The positive constant $\sigma$ is the volatility of $S$. The specific form of the otherwise arbitrary constant $\varpi$ depends on the nature of the security modelled which could be a stock, a currency, a commodity etc. For example, if $S$ is a stock paying a dividend at the continuous rate $\delta$, we have $\varpi = r - \delta$.

**4. Asian options and their equilibrium pricing:** In the Black–Scholes framework of §3, fix any time $t_0$ and consider the accumulation process $J$ given for any time $t$ by:

$$J(t) = \int_{t_0}^{t} S_u \, du\,,$$

The *arithmetic–average Asian option* written at time $t_0$ with maturity $T$ and strike price $K$ is then the stochastic process on the closed time interval from $t_0$ to $T$ paying

$$\left( \frac{J(T)}{T - t_0} - K \right)^+ = \max \left\{ 0, \frac{J(T)}{T - t_0} - K \right\}$$

at time $T$ and paying nothing at all other times. As such it is a contingent claim on the time interval from $t_0$ to $T$ with payoff $(J(T)/(T - t_0) - K)^+$.

It is one of the fundamental insights that in the equilibrium framework of the Black–Scholes model any such contingent claim on a risky security has an equilibrium price too which is equal to the expectation of its payoff with respect to the risk neutral measure conditional on today's information, see [**MR**, Corollary 5.1.1] for instance. Applying this *arbitrage pricing principle*, the price $C_t$ of the Asian option at any time $t$ between $t_0$ and $T$ so is given as the $Q$–expectation conditional on the information $\mathscr{F}_t$ available at time $t$:

$$C_t = e^{-r(T-t)} E^Q \left[ \left( \frac{J(T)}{T - t_0} - K \right)^+ \Big| \mathscr{F}_t \right].$$

However, following [**GY**, §3.2], do not focus on this price, but normalize the valuation problem as follows. On factoring out the reciprocal of the lenght $T - t_0$ of the time period, split the integral $J(T)$ into two integrals, one of which is deterministic by time $t$ and the other of which is random. Couple the deterministic integral with the new strike. For the random integral, restart the Brownian motion driving the underlying at time $t$, and then using the scaling property of Brownian motion, change time to normalize its coefficient in the new time scale to two. The precise result is the factorization:

$$C_t = \frac{e^{-r(T-t)}}{T - t_0} \cdot \frac{4S_t}{\sigma^2} \cdot C^{(\nu)}(h, q)\,,$$

which reduces the general valuation problem to computing

$$C^{(\nu)}(h, q) = E^Q \left[ (A_h^{(\nu)} - q)^+ \right],$$

the *normalized time–t price* of the Asian option. To explain the concepts, $A^{(\nu)}$ is Yor's accumulation process

$$A_h^{(\nu)} = \int_0^h e^{2(B_w+\nu w)} dw \,,$$

and the normalized parameters are as follows:

$$\nu = \frac{2\varpi}{\sigma^2} - 1, \qquad h = \frac{\sigma^2}{4}(T-t), \qquad q = kh + q^*,$$

where

$$k = \frac{K}{S_t}, \qquad q^* = q^*(t) = \frac{\sigma^2}{4S_t}\left(K \cdot (t-t_0) - \int_{t_0}^t S_u \, du\right).$$

To interpret these quantities, $\nu$ is the *normalized adjusted interest rate*, $h$ is the *normalized time to maturity*, which is non–negative, and $q$ is the *normalized strike price*.

# Part II Statement of results

**5. A reduction of the valuation problem:** Computing the normalized time–$t$ price of the Asian option reduces to the case where the normalized strike price $q$ is positive. Indeed, if $q$ is non–positive, Asian options loose their option feature, and their normalized time–$t$ price is given by

$$C^{(\nu)}(h,q) := E^Q\big[\big(A_h^{(\nu)} - q\big)^+\big] = E^Q\big[A_h^{(\nu)}\big] - q \,.$$

On applying Fubini's theorem, this last expectation is computed as follows:

$$E^Q\left[A_h^{(\nu)}\right] = \frac{e^{2h(\nu+1)} - 1}{2(\nu+1)} \,,$$

for $\nu$ any real number, and being equal to $h$ if $\nu = -1$ in particular.

**6. Statement of the main results:** The main result of this paper is a closed form solution for the normalized price $C^{(\nu)}(h,q)$ of the Asian option in the generic case where $q$ is positive. It expresses this function as the sum of integral representations. These are obtained by integrating the product of Hermite functions $H_\mu$, discussed in §2, with functions derived from weighted complementary error functions. Proved later in Part III, the precise result is as follows:

**Theorem:** *If $q$ is positive, the normalized price $C^{(\nu)}(h,q)$ of the Asian option is given by the following difference*

$$C^{(\nu)}(h,q) = ce^{2h(\nu+1)}\, S_{\nu+2} - c\, S_\nu$$

*where the $S_\xi$ are three–term sums*

$$S_\xi = C_{\mathrm{trig},\theta_\xi,\xi}(\rho_\xi) + C_{\mathrm{hyp},\theta_\xi,\xi}(\rho_\xi) + C_{\mathrm{hyp},\theta_\xi,-\xi}(\rho_\xi)$$

*whose single summands are integrals that depend on $\theta_\xi \in [\frac{\pi}{2}, \pi]$ and $\rho_\xi \geq 0$, but which as a whole are independent of these parameters.*

To explain the concepts, $c$ is given by

$$c = c(\nu, q) = \frac{\Gamma(\nu+4) \cdot (2q)^{\frac{1}{2}(\nu+2)}}{2\pi \cdot (\nu+1) \cdot e^{\frac{1}{2q}}}$$

recalling $\nu = 2\varpi/\sigma^2 - 1$. With $\theta$ any real in the closed interval $[\frac{\pi}{2}, \pi]$ and $\rho$ any non–negative real, the *trigonometric terms* $C_{\mathrm{trig},\theta,\xi}(\rho)$ with $\xi$ equal to $\nu$ or $\nu+2$ are the integrals

$$C_{\mathrm{trig},\theta,\xi}(\rho) = \int_0^\theta \mathrm{Re}\left(H_{-(\nu+4)}\left(-\frac{\cosh(\rho+i\phi)}{\sqrt{2q}}\right)E_\xi(h)(\rho+i\phi)\right)d\phi$$

over the real parts of products of Hermite functions times certain functions $E_b(h)$, and the *hyperbolic terms* $C_{\mathrm{hyp},\theta,\xi}(\rho)$ with $\xi$ equal to $\pm\nu$ or $\pm(\nu+2)$ are the integrals

$$C_{\mathrm{hyp},\theta,\xi}(\rho) = \int_\rho^\infty \mathrm{Im}\left(H_{-(\nu+4)}\left(-\frac{\cosh(y+i\theta)}{\sqrt{2q}}\right)E_\xi(h)(y+i\theta)\right)dy$$

over the imaginary parts of such products. Herein $E_\xi(h)$ are the weighted complementary error functions for any complex $w$ given by

$$E_\xi(h)(w) = e^{w\xi}\mathrm{Erfc}\left(\frac{w}{\sqrt{2h}} + \frac{\xi}{2}\sqrt{2h}\right).$$

**Remark:** If $\rho$ equals zero, the trigonometric terms specialize to

$$C_{\mathrm{trig},\theta,\xi}(0) = 2\int_0^\theta H_{-(\nu+4)}\left(-\frac{\cos(\phi)}{\sqrt{2q}}\right)\cos(\xi\phi)\,d\phi\,.$$

Notice that for $\nu+4$ any non–positive integer, the Hermite functions of the formula specialize to the corresponding Hermite polynomials. As explained in [**SH**, §7] and [**ST**,§4 Corollary], the normalized price is then given by a finite number of weighted complementary error functions and their derivatives. This effect is also observed in [**Du**] where Asia densities are expressed as certain integrals against Hermite functions.

Finally deriving benchmarks for the normalized prices of Asian options seems to be one of the main practical application of formulas as that of the Theorem. To get an impression, Asian option on stocks are usually written with maturities of one year and issued at par: $K = S_0$ taking $t_0 = 0$. Current annual volatilities for stocks range between 20 and 50 percent. Taking an annual interest rate of nine percent, we then have from [**ST**, §12]

| $\sigma$ | maximal error | $C^{(\nu)}(h,q)$ |
|---|---|---|
| 20% | $4.9727 \times 10^{-16}$ | 0.00074155998788343 |
| 30% | $4.9687 \times 10^{-16}$ | 0.00217354504625037 |
| 40% | $4.9157 \times 10^{-16}$ | 0.00478100328341654 |
| 50% | $4.9461 \times 10^{-16}$ | 0.00890942045213227 |

**Table 1.** Normalized prices $C^{(\nu)}(h,q)$ of the Asian option for $T{=}1$.

as benchmark normalized prices for the Asian option sharpening the results of [**RS**].

**7. Yor's triple integral and a comparison:** The closed form results of §6 have all been preceded by Yor's integral representation [**Y**, (6.e), p.528]. For the normalized time–$t$ price $C^{(\nu)}(h,q)$ the latter implies the following representation as triple integral:

$$C^{(\nu)}(h,q) = c_{\nu,h} \int_0^\infty x^\nu \int_0^\infty e^{-\frac{(1+x^2)y}{2}} \cdot \left(\frac{1}{y} - q\right)^+ \cdot \psi_{xy}(h)\, dy\, dx\,,$$

setting

$$c_{\nu,h} = \frac{1}{\pi\sqrt{2h}}\, e^{\frac{\pi^2}{2h} - \frac{\nu^2 h}{2}}$$

and where for any positive real number $a$, the function $\psi_a$ is given for any $h > 0$ by

$$\psi_a(h) = \int_0^\infty e^{-\frac{w^2}{2h}} e^{-a\cosh(w)} \sinh(w) \cdot \sin\left(\frac{\pi}{h} w\right) dw\,.$$

While Yor's formula seems to require $\nu$ to be bigger than at least minus one, it is valid for all reals $\nu$. This is proved in [**SL**, §6]. Our formula, in contrast, is given as a sum of single integrals whose integrands have a structural interpretation as products of two functions. It identifies the higher transcendental functions occuring as factors in these products, and shows how they are given by or built up from Hermite functions. On a technical level, these differences can be regarded as consequences of the different mathematical approaches for proving the valuation formula. In fact, Yor's originates from a direct attack on the Asia density whereas our's is eventually based on the indirect enveloping construction of [**CS**], and thus on non–Asian options.

However, there are not only structural differences between the two results. One purpose of explict pricing formulas is to provide means for actually computing option prices in real life. And here Yor's formula has a number of structural difficulties. For instance in the setting of the §6 example we compute the factors $c_{\nu,h}$ as follows:

| $c_{\nu,h}$ | 20% | $\sigma = 30\%$ | $\sigma = 40\%$ | $\sigma = 50\%$ |
|---|---|---|---|---|
| $T = 1$ year | $4.379 \times 10^{214}$ | $2.647 \times 10^{95}$ | $4.266 \times 10^{53}$ | $1.753 \times 10^{34}$ |
| $T = 6$ months | $1.321 \times 10^{429}$ | $6.717 \times 10^{190}$ | $2.288 \times 10^{107}$ | $4.846 \times 10^{68}$ |

**Table 2.** $c_{\nu,h}$ as function of $T$ and $\sigma$.

Yor's formula thus expresses the normalized price of Asian options, which is not too big, as the product of a big number times a triple integral. The latter so has to be small. And it must be computed with very high accuracies to get reasonably accurate results. The formulas of §? Theorem seem to be better behaved in this respect. For instance, in the case $T = 1$ and $\sigma = 30\%$ we have from [**SH**, §11] as values for the hyperbolic terms where $c = 4.072 \times 10^{-12}$. And using any of the formulas of [**ST**, Part II] we get 0.00104 as total contribution of the trigonometric terms.

| $\xi$ | $c \cdot C_{\mathrm{hyp},\frac{\pi}{2},\xi}(0)$ |
|---|---|
| $\nu$ | $-1534170.07497$ |
| $-\nu$ | $1534169.98118$ |
| $\nu+2$ | $-4198799.82516$ |
| $-(\nu+2)$ | $4198799.74048$ |

**Table 3.** Values of $c \cdot C_{\mathrm{hyp},\frac{\pi}{2},\xi}(0)$ for $T{=}1$ and $\sigma{=}30\%$.

# Part III Proof of the valuation formula

**8. Reduction to the Laplace inversion of certain Laplace transforms:** The Laplace transform has been forged by Yor and his school into a very effective tool for transferring stochastics problems into analysis, and [**Y01**] is the last addition to his work on averaging problems. It seems to be characteristic for these problems that it is Bessel function which enter at the Laplace transform level. While this is briefly reviewed too, this section's aim is more precisely to reduce computing the normalized time–$t$ price

$$C^{(\nu)} := E^Q\Big[\big(A_h^{(\nu)} - (kh+q^*)\big)^+\Big]$$

of the Asian option introduced in §4 to the following Laplace inversion problem.

**Proposition:** *If the positive real number $h$ is such that $q(h) = kh+q^*$ is positive, the normalized price of the Asian option at $h$ is given by:*

$$C^{(\nu)} = \mathscr{L}^{-1}\Big(F_{GY}\big(q(h), z\big)\Big)(h)$$

*as the Laplace inverse of $F_{GY}(q(h), z)$ at $h$.*

Here we use the concepts of §4. To explain the functions, for any real $a > 0$ define:

$$F_{GY}(a, z) = \frac{D_\nu(a, z)}{z \cdot \big(z - 2(\nu+1)\big)}\,,$$

for any complex number $z$ with positive real part bigger than $2(\nu+1)$, where on choosing the principal branch of the logarithm:

$$D_\nu(a, z) = \frac{e^{-\frac{1}{2a}}}{a} \int_0^\infty e^{-\frac{x^2}{2a}} \cdot x^{\nu+3} \cdot I_{\sqrt{2z+\nu^2}}\left(\frac{x}{a}\right) dx\,.$$

Here $I_\mu$ is, for any complex number $\mu$, the modified Bessel function of order $\mu$ more fully discussed in §7 below or in [**L**, Chapter 5], for instance.

The basic idea originating with Yor for proving such results is to make the time $h$ a variable and compute the Laplace transform of the functions so obtained. In the case of $C^{(\nu)}$ this is made difficult by the time variable entering not only via the stochastic process but also via the normalized strike price. Following [**C**], as a first step in proving the Proposition so consider for any real number $a$ the non–Asian option price functions $f_{GY,a}$ on the positive real line that send any $x > 0$ to

$$f_{GY,a}(x) = E^Q\big[(A_x^{(\nu)} - a)^+\big].$$

Taken individually they cannot be used to value the original Asian option. However, as a whole they allow one to recover $C^{(\nu)}$. Indeed, if $q = kh+q^*$ is positive, $C^{(\nu)}$ is obtained by choosing the function $f_{GY,kh+q^*}$ and evaluating it at $h$. Equivalently, using the injectivity of the Laplace transform on continuous functions, $C^{(\nu)}$ so is equal to the Laplace inverse at $h$ of the Laplace transform of $f_{GY,kh+q^*}$. Thus we are reduced to identify any $F_{GY,a}$ as the Laplace transform $\mathscr{L}(f_{GY,a})$ of $f_{GY,a}$.

For explaining how Bessel functions enter into the picture, we restrict to the case $\nu \geq 0$ where $\mathscr{L}(f_{GY,a})$ has been computed in [**GY**, §3]. This will illustrate the essential stochastic ideas, and is also sufficient for our purposes. A proof for arbitrary $\nu$ using complex analytic tools in particular can be found in [**CS**]. The punchline of Yor's computations is as follows. The basic idea is to make time stochastic using the *Lamperti identity*

$$e^{\nu w + B_w} = R^{(\nu)}\big(A^{(\nu)}_w\big)$$

for $\nu \geq 0$ proved in [**Y92a**, §2]. Herein $R^{(\nu)}$ is the *Bessel process* of index $\nu$ and with $R^{(\nu)}(0) = 1$, see [**RY**, XI §1]. Then, $A^{(\nu)}$ has the double role of both underlying and stochastic clock. Using the fact that $a$ is independent of time, thus transcribe the condition on the underlying to be bigger than $a$ as the first passage time

$$\tau_{\nu,a} = \inf\{w | A^{(\nu)}_w > a\}$$

for the stochastic clock. This is the key idea for obtaining the representation

$$f_{GY,a}(w) = E^Q\left[\frac{e^{2(\nu+1)[w-\tau_{\nu,a}]^+} - 1}{2(\nu+1)} \cdot (R^{(\nu)}_a)^2\right],$$

for all $w > 0$. At first sight this may have complicated the problem. However, it is just what fits the Laplace transform $L$ of $f_{GY,a}$:

$$L(z) := \int_0^\infty e^{-zw} E^Q\left[\frac{e^{2(\nu+1)[w-\tau_{\nu,a}]^+} - 1}{2(\nu+1)} \cdot (R^{(\nu)}_a)^2\right] dw\,.$$

For computing this integral to show $L = F_{GY,a}$ it seems best to follow a communication of Yor's, switch to measurable functions and applying Tonelli's theorem interchange the Laplace integral with the expectation $E^Q$. In this way we get the identity

$$L(z) = \frac{1}{z(z-2(\nu+1))} E^Q\left[e^{-z\tau_{\nu,a}} (R^{(\nu)})^2\right]$$

of measurable functions for any complex number $z$ with $\mathrm{Re}\,(z) > 2(\nu+1)$. To identify the expectation in the numerator as $D_{\nu,a}$ write it out as

$$\int_0^\infty x^2 \cdot E^Q\left[e^{-z\tau_{\nu,a}}\middle| R^{(\nu)}_a = x\right] \cdot p_{\nu,a}(x)\, dx\,,$$

where $p_{\nu,a}$ is the Bessel semigroup of index $\nu$ starting at 1 at time $a$. Explicit expressions for them are classically known for $\nu > -1$, see [**Y80**, (4.3), p.78] or [**GY**, Proposition 2.2], while the case $\nu < -1$ has been addressed only recently in [**YGö**, §3]. Under the hypothesis $\nu \geq 0$, Yor has in [**Y80**, Théorème 4.7, p.80] (see also [**GY**, Lemma 2.1 and Proposition 2.6]) explicitly computed the conditional expectation factor of this integrand while for $\nu < 0$ there are no results. The upshot is that these results are available both only if $\nu \geq 0$, and this is where the restriction on $\nu$ comes from. On substitution we get

$$D_\nu(a,z) = \frac{e^{-\frac{1}{2a}}}{a} \int_0^\infty e^{-\frac{x^2}{2a}} x^{\nu+3} I_{\sqrt{2z+\nu^2}}\left(\frac{x}{a}\right) dx\,,$$

as required. To complete the Tonelli argument proposed by Yor and to complete the proof, we have to establish the finiteness of this integral for any fixed complex number $z$ with $\mathrm{Re}\,(z) > 2(\nu+1)$. There is a further technical point to be taken care of. Choose the principal branch of the logarithm to define the square root on the complex plane with the non–positive real line deleted, whence $\mu = (2z+\nu^2)^{1/2}$ has a positive real part. Integrability now follows using the asymptotic behaviour of the Bessel function $I_\mu$ near the origin and towards infinity, and we have shown $L = F_{GY,a}$ if $\nu \geq 0$ as required.

**9. Preliminaries on integral representations of gamma and Bessel functions:** With hindsight, it is the gamma function which is at the base of valuing the Asian option. However, the gamma function does not enter via Euler's integral recalled to be given by

$$\Gamma(s) = \int_0^\infty e^{-x} x^{s-1}\, dx\,,$$

for any complex $s$ with $\mathrm{Re}\,(s) > 0$. The problem here is that such a formula would be required for $s$ with non–positive real parts where the resulting integrand is not integrable near the orgin. Being able to perform such integrations from minus infinity to plus infinity avoiding the origin needs higher dimensions, and recall Hankel found natural candidates for this some 130 years after Euler. His *Hankel contours* $C_{\theta,R}$ with parameters $\theta$ any angle in $(\frac{\pi}{2}, \pi)$ and $R > 0$ any real are the following contours of integration in the complex plane. One comes in from infinity on the ray $\{y \exp(-i\theta) : y \geq R\}$ until the point $R \exp(i\theta)$, then passes counterclockwise around zero with distance $R$ until the point $R\exp(i\theta)$, and leaves on the ray $\{y \exp(i\theta) : y \geq R\}$ to infinity. And using integration over any of these contours, we then have *Hankel's formula*:

$$\frac{1}{\Gamma(s)} = \frac{1}{2\pi i} \int_C e^{\xi} \cdot \xi^{-s}\, d\xi\,,$$

valid now for any complex number $s$, see for instance [**D**, pp.225f].

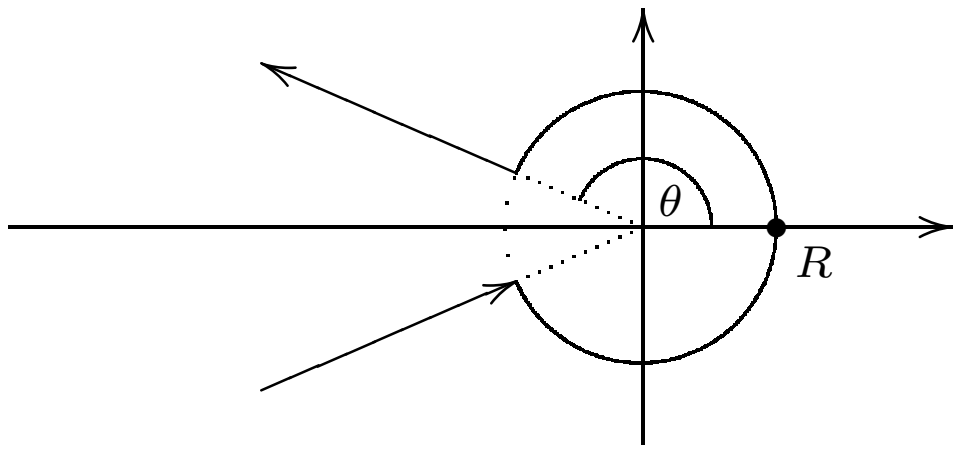

Figure 1. The Hankel contour $C_{\theta,R}$.

The occurrence of modified Bessel functions in §6 Proposition now furnishes the point of attack. From [**L**, §5.7] they are given for any complex number $\mu$ by the following series:

$$I_\mu(z) = \sum\nolimits_{m=0}^{\infty} \frac{1}{m!\,\Gamma(\mu+m+1)} \cdot \left(\frac{z}{2}\right)^{\mu+2m}\,,$$

for any complex $z$ in $\mathbf{C}\backslash\mathbf{R}_{<0}$. The idea is to substitute for any of the reciprocal gamma factors of the single terms of this seris substitute the corresponding Hankel formulas and interchange the order of summation and contour integration. With the details in [**WW**, 17·231, p.362], for instance, the result is the following Hankel–type integral representations of modified Bessel functions on the right half–plane

**Lemma:** *For any modified Bessel function $I_\mu$ one has:*

$$I_\mu(z) = \frac{1}{2\pi i} \int_{\log C} e^{-\mu\cdot\xi + z\cdot\cosh(\xi)}\, d\xi\,,$$

*for any complex number $z$ with positive real part and any Hankel contour $C$.*

Herein the principal branch of the logarithm on the complex plane with the non–positive real axis deleted has been chosen. The contour $\log C_{\theta,R}$ thus has the following shape: Coming in from plus infinity, move on the parallel through the point $-i\theta$ to the real axis to the point $\log R - i\theta$. From this point move up to the point $\log R + i\theta$ on a parallel to the imaginary axis. Finally exit from $\log R + i\theta$ to plus infinity on the parallel through $i\theta$ to the real axis. Notice $\log C_R \subseteq \{z|\mathrm{Re}\,(z) > 0\}$ if and only if the radius of the circle in $C_{\theta,R}$ is bigger than 1.

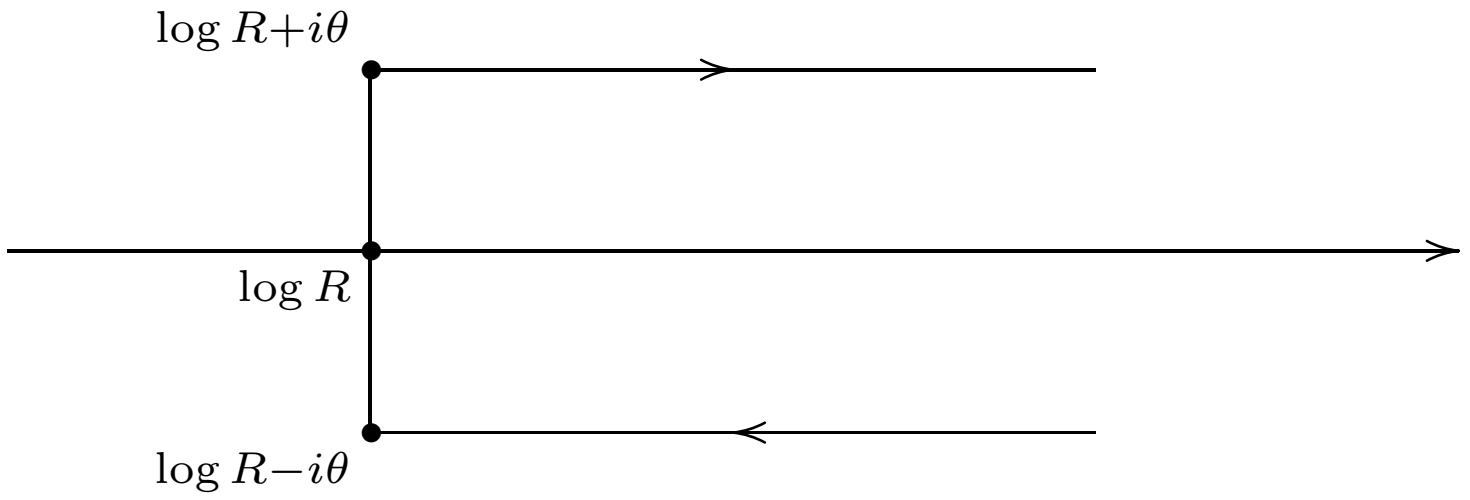

Figure 2. The contour $\log C_{\theta,R}$.

**10. First steps of the Laplace inversion:** This section resumes the proof of §6 Theorem by addressing the Laplace inversion problem of §8 Proposition. Define for any positive real $a$ and for $b$ equal to $\nu$ or $\nu+2$ functions on the positive real line by sending any $h$ to

$$\Sigma^a_b(h) = \int_0^\infty \frac{1}{2\pi i} \int_{\log C_{\theta,R}} e^{-x^2 + x\frac{\sqrt{2}}{\sqrt{a}}\cosh(w)}\, x^{\nu+3}\, \mathscr{L}^{-1}\left(\frac{e^{-w\sqrt{z}}}{z-(\nu+2)^2}\right) dw\, dx\,,$$

where $\log C_{\theta,R}$ is any logarithmicalized Hankel contour as in §7, and let $c_1$ denote

$$c_1 = e^{-\frac{1}{2}\nu^2 h} \cdot \frac{(2a)^{\frac{1}{2}(\nu+2)}}{(\nu+1)\cdot e^{\frac{1}{2a}}}\,.$$

Then §8 Proposition reduces to compute the functions $\Sigma^a_b$ as follows

**Lemma:** *If $\nu$ is non–negative, the inverse Laplace transform at any positive real number $h$ of any function $F_{GY}(a,z)$ is given by:*

$$\mathscr{L}^{-1}\Big(F_{GY}(a,\cdot)\Big)(h) = c_1\Big(\Sigma^a_{\nu+2}(h) - \Sigma^a_\nu(h)\Big)\,,$$

*where any of the functions $\Sigma^a_b$ is independent of the Hankel contour chosen and has an absolutely integrable integrand.*

The proof of the Lemma is based on the complex inversion formula for the Laplace transform, see [**B**, Chapter 7, §5].

To prove the Lemma, first recall that for any positive real number $a$ we have:

$$F_{GY}(a,z) = \frac{D_\nu(a,z)}{z \cdot \big(z - 2(\nu+1)\big)} ,$$

where

$$D_\nu(a,z) = \frac{e^{-\frac{1}{2a}}}{a} \int_0^\infty e^{-\frac{x^2}{2a}} \cdot x^{\nu+3} \cdot I_{\sqrt{2z+\nu^2}}\left(\frac{x}{a}\right) dx ,$$

for any complex number $z$ with real part bigger than $2(\nu{+}1)$. Then let $z_0'$ be any positive real number such that the line $\{z|\mathrm{Re}\,(z) = z_0'\}$ is contained in the half–plane where $F_{GY}(a,z)$ is a holomorphic function. Suppose for a moment proved that this function is rapidly decreasing, and apply the complex inversion formula to it. Writing out $D_\nu(a,z)$, the inverse Laplace transform of $F_{GY}(a,z)$ is the function on the positive real line given for any positive real number $h$ by:

$$\frac{e^{-\frac{1}{2a}}}{a} \cdot \frac{1}{2\pi i} \int_{z_0'-i\infty}^{z_0'+i\infty} e^{hz} \int_0^\infty \frac{e^{-\frac{x^2}{2a}} \cdot x^{\nu+3} \cdot I_{\sqrt{2z+\nu^2}}\left(\frac{x}{a}\right)}{z(z-2(\nu{+}1))} \, dx \, dz .$$

For an equivalent expression, change variables $\eta = 2z{+}\nu^2$, put $z_0 = 2z_0'{+}\nu^2$, and substitute the Hankel–type integral representation of §9 Lemma for the modified Bessel function. With $c_1' = 4(\nu+1) \cdot (2a)^{-(\nu+4)/2} \cdot c_1$, the following integral is then to be computed:

$$c_1' \left(\frac{1}{2\pi i}\right)^2 \int_{z_0-i\infty}^{z_0+i\infty} \int_0^\infty \int_{\log C} e^{-\frac{x^2}{2a}} x^{\nu+3} e^{\frac{x}{a}\cosh(w)} \cdot e^{\frac{h}{2}\cdot z} \frac{e^{-w\sqrt{z}}}{(z-\nu^2)(z-(\nu{+}2)^2)} \, dw \, dx \, dz .$$

We claim that for $R$ sufficiently big the absolute value of the integrand of this integral is exponentially decreasing to zero with the absolute values of $z$, $x$, or $w$ going to infinity.

Granting this result, the above triple integral then gives the desired Laplace inverse, and, using Fubini's theorem, the order of its integrals can be interchanged. Take the integral for Laplace inversion, i.e., the integral over the line $\{z|\mathrm{Re}\,(z) = z_0\}$, as inner integral. Change variables $x = (2a)^{1/2}t$. The Lemma follows on decomposing the denominator of the integrand in partial fractions.

One is thus reduced to prove the last claim about the asymptotic behaviour. For the calculations recall $|\exp(\xi)| = \exp(\mathrm{Re}\,(\xi))$, for any complex number $\xi$. For the asymptotic behaviour in the absolute value of $w$, reduce to elements $w = x \pm i\theta$ in $\log C_{\theta,R}$. For these $\mathrm{Re}\,(\cosh(w)) = -\cosh(x)|\cos(\theta)|$. The absolute value of the hyperbolic cosine factor of the numerator equals $\exp((x/a)\mathrm{Re}\,(\cosh(w))$, whence the required asymptotic behaviour in the absolute value of $w$.

For the asymptotic behaviour in $z$, notice $|\exp(-wz^{1/2})| = \exp(-\mathrm{Re}(wz^{1/2}))$ and recall $z^{1/2} = \exp((1/2)\cdot(\log|z| + i\arg(z)))$. The argument of $z$ converges to $\pi/2$ with $|z|$ going to infinity. Now $\mathrm{Re}(wz^{1/2})$ equals $|z|^{1/2}$ times $\mathrm{Re}(w)\cos(\arg(z)/2) - \mathrm{Im}(w)\sin(\arg(z)/2)$. Herein the cosine is positive and bigger than $\cos(\pi/4)$. Thus $\mathrm{Re}(wz^{1/2})$ is positive if the real part of any $w$ is big enough. Hence choose $R$ big enough to have the desired asymptotic behaviour in the absolute value of $z$.

The asymptotic behaviour in $x$ is determined by $\exp(-x^2/2a)$, thus completing the proof.

**11. Computation of certain Laplace transforms:** This section computes the inverse Laplace transforms identified in §10 Lemma thereby illustrating typical techniques for working with Laplace transforms.

For any $\alpha$ and $\beta$ in $\mathbf{C}$, consider the functions on the positive real line given by:

$$f_{\alpha,\beta}(t) = \frac{t^{-1/2}}{\sqrt{\pi}} e^{-\frac{\alpha^2}{4\cdot t}} - \beta\cdot e^{\alpha\cdot\beta+\beta^2\cdot t}\,\mathrm{Erfc}\left(\beta\sqrt{t} + \frac{\alpha}{2\sqrt{t}}\right),$$

$$g_{\alpha,\beta}(t) = \frac{e^{\beta^2\cdot t}}{2}\left(e^{+\alpha\cdot\beta}\,\mathrm{Erfc}\left(\frac{\alpha}{2\sqrt{t}} + \beta\sqrt{t}\right) + e^{-\alpha\cdot\beta}\,\mathrm{Erfc}\left(\frac{\alpha}{2\sqrt{t}} - \beta\sqrt{t}\right)\right),$$

for any positive real number $t$. Then one has the following two results:

**Lemma:** *If the real parts of $\alpha$ and $\alpha^2$ are positive, one has:*

$$\mathscr{L}(f_{\alpha,\beta})(z) = \frac{e^{-\alpha\sqrt{z}}}{\sqrt{z} + \beta},$$

*for any complex number $z$ in $\mathbf{C}\setminus\mathbf{R}_{<0}$ with real part bigger than* $|\mathrm{Re}(\beta)|$.

**Corollary:** *If the real parts of $\alpha$ and $\alpha^2$ are positive, one has:*

$$\mathscr{L}(g_{\alpha,\beta})(z) = \frac{e^{-\alpha\sqrt{z}}}{z - \beta^2},$$

*for any complex number $z$ in $\mathbf{C}\setminus\mathbf{R}_{<0}$ with real part bigger than* $\mathrm{Re}(\beta^2)$.

The Corollary follows from the Lemma upon decomposing the denominator in partial fractions and using the linearity of the Laplace transform.

We use the following two results proved mutatis mutandis in [**D**, Beispiel 8, p.50f]:

$$\mathscr{L}(\psi_\alpha)(z) = e^{-\alpha\sqrt{z}} \qquad \text{where} \qquad \psi_\alpha(t) = \frac{\alpha}{2\sqrt{\pi}}\cdot t^{-\frac{3}{2}}\cdot e^{-\frac{\alpha^2}{4\cdot t}},$$

$$\mathscr{L}(\chi_\alpha)(z) = \frac{e^{-\alpha\sqrt{z}}}{\sqrt{z}} \qquad \text{where} \qquad \chi_\alpha(t) = \frac{1}{\sqrt{\pi}}\cdot t^{-\frac{1}{2}}\cdot e^{-\frac{\alpha^2}{4\cdot t}},$$

for any complex number $z$ in $\mathbf{C}\setminus\mathbf{R}_{\leq 0}$ and any positive real number $t$.

Subtracting $\chi_\alpha$ from $f_{\alpha,\beta}$, the proof of the Lemma reduces to show the identity:

$$\mathscr{L}\left(-\beta\cdot e^{\alpha\cdot\beta+\beta^2\cdot t}\,\mathrm{Erfc}\left(\beta\sqrt{t}+\frac{\alpha}{2\sqrt{t}}\right)\right)(z)=-\beta\frac{e^{-\alpha\sqrt{z}}}{\sqrt{z}(\sqrt{z}+\beta)}.$$

Multiplying any nice function $f$ with $\exp(a\cdot)$ induces a shift by $-a$ in its Laplace transform: $\mathscr{L}(\exp(at)f(t))(z)=\mathscr{L}(f)(z-a)$. Using this with $a=\beta^2$ one is further reduced to calculating the Laplace transform of the above complementary error function factor only.

Since $\mathrm{Re}(\alpha)$ is positive, the real part of $\beta u^{1/2}+(\alpha/2)u^{-1/2}$ goes to plus infinity with $u$ converging from the right to zero. The Fundamental Theorem of Calculus thus gives:

$$\mathrm{Erfc}\left(\beta\sqrt{t}+\frac{\alpha}{2\sqrt{t}}\right)=-\frac{2}{\sqrt{\pi}}\int_0^t e^{-\left(\beta\sqrt{u}+\frac{\alpha}{2\sqrt{u}}\right)^2}\left(\frac{\beta}{2}u^{-1/2}-\frac{\alpha}{4}u^{-3/2}\right)du\,.$$

Using the transform–of–an–integral property $\mathscr{L}(\int_0^\bullet f(u)\,du)(w)=w^{-1}\mathscr{L}(f)(w)$ the Laplace transform of this complementary error function at $w=z-\beta^2$ is given by:

$$\begin{aligned}&-\frac{2\beta}{2(z-\beta^2)}\mathscr{L}\left(\frac{1}{\sqrt{\pi}}\cdot t^{-\frac{1}{2}}\cdot e^{-\left(\beta\sqrt{u}+\frac{\alpha}{2\sqrt{u}}\right)^2}\right)(z-\beta^2)\\&+\frac{2}{2(z-\beta^2)}\mathscr{L}\left(\frac{\alpha}{2\sqrt{\pi}}\cdot t^{-\frac{3}{2}}\cdot e^{-\left(\beta\sqrt{u}+\frac{\alpha}{2\sqrt{u}}\right)^2}\right)(z-\beta^2).\end{aligned}$$

Using the Laplace transform of $\chi_\alpha$, the first Laplace transform of this sum equals

$$\frac{e^{-\alpha\beta}e^{-\alpha\sqrt{z}}}{\sqrt{z}}.$$

Using the Laplace transform of $\psi_\alpha$, the second Laplace transform of this sum equals

$$e^{-\alpha\beta}e^{-\alpha\sqrt{z}}.$$

The above sum hence equals $\exp(-\alpha\beta)\exp(-\alpha\sqrt{z})/(z^{1/2}(z^{1/2}+\beta))$. The identity to be proved follows upon substituting this last expression. This completes the proof.

**12. Two intermediate results:** In this section the Laplace inversion of §10 is resumed concentrating on the two summands $\Sigma_b^a$ of the expression derived in §10 Lemma. If $\nu$ is non–negative, recall them as the functions on the positive real line sending any $h$ to

$$\Sigma_b^a(h)=\frac{1}{2\pi i}\int_{\log C_{\theta,R}}\int_0^\infty e^{-x^2+x\frac{\sqrt{2}}{\sqrt{a}}\cosh(w)}x^{\nu+3}\mathscr{L}^{-1}\left(\frac{e^{-w\sqrt{z}}}{z-b^2}\right)\left(\frac{h}{2}\right)dx\,dw$$

where $a>0$ and $b$ are any reals, and $\log C_{\theta,R}$ is any logarithmicalized Hankel contour as discussed in §9. For these integrals we have the following key result:

**Proposition:** *If moreover $R \geq 1$, we have*

$$\Sigma_b^a(h) = \frac{1}{2}\Gamma(\nu+4)e^{\frac{1}{2}b^2h} \cdot \frac{1}{2\pi i}\int_{\log C_{\theta,R}} \Big(E_b(h) + E_{-b}(h)\Big)(w) \cdot H_{-(\nu+4)}\Big(-\frac{\cosh(w)}{\sqrt{2a}}\Big)\, dw$$

*where the respective integrands on both sides are absolutely integrable, and the integral on the right hand side does not depend on the Hankel contour chosen.*

Notice that the $\Sigma_b^a$ are complex–valued integrals, and we so have the decomposition

$$\Sigma_b^a(h) = \mathrm{Re}\,(\Sigma_b^a)(h) + i\,\mathrm{Im}\,(\Sigma_b^a)(h)$$

into real and imaginary parts. Our second result then expresses the contour integral representation of the Proposition in terms of these real–valued integrals as follows

**Corollary:** *Under the additional assumption $\rho = \log R \geq 0$, we have*

$$\mathrm{Re}\,(\Sigma_b^a)(h) = \frac{1}{2\pi}e^{\frac{1}{2}b^2h}\int_0^\infty e^{-x^2}x^{\nu+3}\int_0^\theta \mathrm{Re}\,\big(F_{x,b}^a(h)\big)(\rho+i\phi)\,d\phi\,dx$$

$$\mathrm{Im}\,(\Sigma_b^a)(h) = \frac{1}{2\pi}e^{\frac{1}{2}b^2h}\int_0^\infty e^{-x^2}x^{\nu+3}\int_\rho^\infty \mathrm{Im}\,\big(F_{x,b}^a(h)\big)(y+i\theta)\,dy\,dx\,.$$

To explain the notation in these two results, the Hermite functions $H_\mu$ of degree $\mu$ are discussed in §2, and the two other functions that occur are given by:

$$E_b(h)(w) = e^{wb}\cdot \mathrm{Erfc}\left(\frac{w}{\sqrt{2h}} + \frac{b}{2}\sqrt{2h}\right),$$

$$F_{x,b}^a(h)(w) = e^{x\frac{\sqrt{2}}{\sqrt{a}}\cosh(w)}\Big(E_b(h)(w) + E_{-b}(h)(w)\Big),$$

for any complex number $w$.

**Proof of the Proposition:** For the proof of the Proposition choose a Hankel contour $C_{\theta,R}$ with $R$ so big that for any element $w$ in $\log C_{\theta,R}$ also the real parts of both $w$ and $w^2$ are positive. At any $w$ in $\log C_{\theta,R}$ substitute in $\mathscr{I}_b^a(h)$ the inverse Laplace transforms of §11 Corollary. Interchange the order of integration using the absolute integrability of the integrand. This gives the expression of the Lemma for $\mathscr{I}_b^a(h)$. The integrand herein is a holomorphic function on $\mathbf{C}\setminus\mathbf{R}_{<0}$. Using the Cauchy Theorem, the value of the integral is independent of the Hankel contour $C_{\theta,R}$ chosen as long as $\log R$ is non–negative. This completes the proof.

**Proof of the Corollary:** Given the Proposition, the proof of the Corollary is an exercise in path integration. Put $F_x = F_{x,b}^a(h)$. Changing the order of integration there,

$$\int_0^\infty e^{-x^2}x^{\nu+3}\frac{1}{2\pi i}\int_{\log C_{\theta,R}} F_x(w)\,dw\,dx$$

is to be calculated. Concentrate on its inner integral and consider the following subpath

$$P = P_c + P_\infty \qquad \text{with} \qquad P_c = P_{c,+} - P_{c,-}$$

of the path $\log C_{\theta,R}$. Here the path $P_{c,-}$ starts from $\log R$ and moves parallel to the imaginary axis to the point $\log R - i\theta$. The path $P_{c,+}$ starts from $\log R$ and moves parallel to the imaginary axis to the point $\log R + i\theta$. The path $P_\infty$ moves from $\log R + i\theta$ parallel to the real axis to $+\infty$. For any $x > 0$, the inner integral thus breaks up as follows:

$$\frac{1}{2\pi i}\int_{P_\infty}\Big(F_x(w) - F_x(\overline{w})\Big)\,dw + \frac{1}{2\pi i}\int_{P_{c,+}} F_x(w)\,dw \; - \frac{1}{2\pi i}\int_{P_{c,-}} F_x(w)\,dw\,,$$

with $\overline{w}$ the complex conjugate of $w$. Using the series expansion of the exponential and the complementary error functions, $F_x$ is compatible with complex conjugation, i.e., $F_x$ evaluated at the complex conjugate of any complex number $w$ is the complex conjugate of $F_x$ at $w$:

$$F_x\left(\overline{w}\right) = \overline{F_x(w)}.$$

Write the elements of $P_\infty$ as $w = y + i\theta$ with $y \geq \log R$, notice $dw = dy$ and change variables. Using the compatibility of $F_x$ with complex conjugation, one obtains the improper integral from $\log R$ to infinity of the imaginary parts of $F_x(y+i\theta)$.

On the circle part $P_c$ of $P$ one has to be a bit more careful about the volume forms. The elements of $P_{c,+}$ are parametrized by $w = \log R + i\phi$ with $\phi$ in $[0, \theta]$, whereas those of $P_{c,-}$ are parametrized by $\log R - i\phi$ with $\phi$ in $[0, \theta]$. Changing variables accordingly and integrating from zero to $\theta$, the induced volume form in the $P_{c,+}$–integral is $dw = i\,d\phi$, whereas that on the $P_{C,-}$–part is $dw = -i\,d\phi$. Using the compatibility of $F_x$ with complex conjugation, this completes the proof of the Corollary.

**13. Explicit calculations:** This section explicitly computes the integrals of §12 Corollary in the case where the parameters $\rho$ there are equal to zero. Recalling for this the function $\Sigma_b^a$ given on the positive real line by

$$\Sigma_b^a(h) = \frac{1}{2\pi i}\int_{\log C}\int_0^\infty e^{-x^2 + x\frac{\sqrt{2}}{\sqrt{a}}\cosh(w)}\, x^{\nu+3}\, \mathscr{L}^{-1}\left(\frac{e^{-w\sqrt{z}}}{z - b^2}\right)\left(\frac{h}{2}\right)\, dx\, dw$$

where $a > 0$ and $b$ are any reals and $C = C_{\theta,1}$ is any Hankel contour as discussed in §9, the precise result to be proved for establishing §6 Remark for $\nu$ non–negative is the

**Lemma:** *If $\nu$ is non–negative, we have*

$$\operatorname{Re}\left(\Sigma_b^a\right)(h) = 2\frac{\Gamma(\nu+4)}{2\pi} e^{\frac{1}{2}b^2 h}\int_0^\theta H_{-(\nu+4)}\left(-\frac{\cos(\phi)}{\sqrt{2a}}\right)\cos(b\phi)\,d\phi\,.$$

The Lemma is proved by computing the integrals over $F_x := F^a_{x,b}$ of §12 Corollary upon choosing $R = 1$ there. For completeness sake start with the imaginary part integrals:

$$\int_0^\infty \mathrm{Im}\,\big(F_x\big)(y+i\theta)\,dy\,.$$

Abbreviate $\beta_\pm = y/\sqrt{2h} \pm (b/2)\sqrt{2h}$, and write the complementary error functions occuring in $F_x(y+i\theta)$ as improper integrals starting from zero as in §2. Multiply out the expressions of the exponents to obtain:

$$e^{\pm(y+i\theta)b}\,\mathrm{Erfc}\Big(\beta_\pm + i\frac{\theta}{\sqrt{2h}}\Big)$$
$$= \frac{2}{\sqrt{\pi}}\cdot e^{\frac{\theta^2}{2h}}\cdot e^{\pm yb}\int_0^\infty e^{-(u+\beta_\pm)^2}\cdot e^{i\theta\Big(\pm b-(u+\beta_\pm)\frac{\sqrt{2h}}{h}\Big)}\,du\,.$$

Since $\beta_\pm$ is a real number and $u$ can be taken as real numbers, the real respectively imaginary parts in $F_x$ are determined by the real respectively imaginary parts of the exponentials in the integral, whence in particular

$$\mathrm{Re}\,\big(E_b(h)\big)(y+i\theta) = \frac{2}{\sqrt{\pi}}e^{\frac{\theta^2}{2h}}\int_{\frac{y}{\sqrt{2h}}+\frac{b}{2}\sqrt{2h}}^\infty e^{-u^2}\cos\Bigg(\pi\Bigg(b-u\frac{\sqrt{2h}}{h}\Bigg)\Bigg)du\,,$$
$$\mathrm{Im}\,\big(E_b(h)\big)(y+i\theta) = \frac{2}{\sqrt{\pi}}e^{\frac{\theta^2}{2h}}\int_{\frac{y}{\sqrt{2h}}+\frac{b}{2}\sqrt{2h}}^\infty e^{-u^2}\sin\Bigg(\pi\Bigg(b-u\frac{\sqrt{2h}}{h}\Bigg)\Bigg)du\,.$$

For proving the Lemma,

$$\int_0^\theta \mathrm{Re}\,\big(F_x\big)(i\phi)\,d\phi$$

is calculated. In contrast to the above argument, in the complementary error functions occuring in $F_x(i\phi)$ now the following paths of integration are used: Abbreviating $\beta_\pm = \pm(b/2)\sqrt{2h}$, first move from $\beta_\pm + i\phi/\sqrt{2h}$ to $\beta_\pm$, then continue from $\beta_\pm$ to plus infinity along the real line. Since $\cosh(i\phi) = \cos(\phi)$, the above integral then equals:

$$\int_0^\pi e^{x\frac{\sqrt{2}}{\sqrt{q}}\cdot\cos(\phi)}\Big\{\big(\mathrm{Erfc}\,(\beta_+) + \mathrm{Erfc}\,(\beta_-)\big)\cdot\cos(\phi b) + \frac{2}{\sqrt{\pi}}\big(\Phi_+ + \Phi_-\big)(\theta)\Big\}\,d\phi,$$

upon abbreviating for any angle $\phi$:

$$\Phi_\pm(\phi) = -\mathrm{Re}\left(e^{\pm i\phi b}\int_{\beta_\pm}^{\beta_\pm+\frac{i\phi}{\sqrt{2h}}} e^{-u^2}\,du\right).$$

In the first of these two last integrals notice $\mathrm{Erfc}\,(\beta_+) + \mathrm{Erfc}\,(\beta_-) = 2$ since $\beta_-$ is minus $\beta_+$. To calculate $\Phi_\pm$ change variables $u = \beta_\pm + iw/\sqrt{2h}$ in the integral. Write the factor $i/\sqrt{2h}$ that is picked up as $\exp(i(\pi/2))/\sqrt{2h}$. Multiplying out the expression obtained in the exponent, it follows

$$e^{\pm i\phi b}\int_{\beta_\pm}^{\beta_\pm+\frac{i\phi}{\sqrt{2h}}} e^{-u^2}\,du \;=\; \frac{e^{-\frac{1}{2}b^2h}}{\sqrt{2h}}\int_0^\phi e^{\frac{u^2}{2h}}\cdot e^{i\left(\frac{\pi}{2}\pm b(\phi-u)\right)}du\,.$$

The real part of this expression is determined by the real part of the exponential functions in the integral. Abbreviating $x_u = b(\theta-u)$, the values of the cosine at $\pi/2\pm x_u$ thus appear as factors. A shift by $\pi/2$ turns a cosine into a sine as follows: $\cos(\pi/2\pm x_u) = \mp\sin(x_u)$. Thus $\Phi_-$ is minus $\Phi_+$, and the expression of the Lemma results.

**14. First part of the proof of the valuation formula:** The proof of the valuation formula of §6 Theorem is in two steps. As a first step, §6 Theorem is in this section established for $\nu$ non–negative. As a second step, this equality is extended in the next section to any complex number $\nu$ using analytic continuation.

Thus let $\nu$ be non–negative. Recalling §10 Lemma, the valuation formula of §6 Theorem is obtained by subtracting $\Sigma^q_\nu(h)$ from $\Sigma^q_{\nu+2}(h)$ and multiplying this difference with the constant $c_1$. Substitute the expressions computed in §12 Proposition. Defining $S_b$ by

$$\Sigma^q_b(h) = \frac{\Gamma(\nu+4)}{2\pi}e^{\frac{1}{2}b^2h}S_b$$

and cancelling the factors $\exp(h\nu^2/2)$ thus completes the proof of the Theorem's representation for the normalized time–$t$ price of the Asia option.

**15. Second part of the proof of the valuation formula:** This second part of the proof of §6 Theorem extends its validity from $\nu$ non–negative, as established in the previous section, to $\nu$ any complex number.

This reduces to show the following two results. First, the normalized price is an entire function in $\nu$. Second, the right hand side of §6 Theorem is a meromorphic function on the complex plane. Indeed, these two functions agree on non–negative real numbers $\nu$. Using the identity theorem they so agree on the complex plane as meromorphic functions. With one of them entire, the other one is then entire, too. In fact, for establishing the second result it is sufficient to divide off the function $c$ with its at most simple poles in the negative integers and prove the following result: The functions $S_b$ of the right hand side of §6 Theorem are meromorphic with at most simple poles off the negative integers.

For proving the first of these two results recall the normalized price $C^{(\nu)}(h,q)$ as given by

$$C^{(\nu)}(h,q) = E\Big[f\Big(A^{(\nu)}(h)\Big)\Big]$$

where $A^{(\nu)}(h) = \int_0^h \exp((2(B_u+\nu u))\,du$, and we have $f(x) = (x-q)^+$, for any real $x$. The proof of this normalized price being entire in $\nu$ is then further reduced to show

$$E\Big[f(A_h)e^{\nu B_h}\Big]$$

an entire function in $\nu$ where $A$ is the process $A^{(0)}$. Indeed, this follows using the Girsanov identity $E[f(A^{(\nu)}(h))] = \exp(-\nu^2h/2)\,E[f(A(h))\exp(\nu B(u))]$ of [**Y**, (1.c), p.510].

For proving the above function entire develop the factor $\exp(\nu B_h)$ of $f(A_h)\exp(\nu B_h)$ into its exponential series. Suppose computing the expectation of the resulting series term by term is justified for any complex number $\nu$. For any complex number $\nu$, this then gives a power series in $\nu$ and thus explicitly shows $E[f(A_h)\exp(\nu B_h)]$ holomorphic at $\nu$.

Interchanging the order of integration and summation is justified using Lebesgue Dominated Convergence if the following is true. The expectations of the absolute value of any single term of the series for $f(A_h)\exp(\nu B_h)$ exist and the series so obtained converges. Using the Cauchy–Schwarz inequality this is implied by the series:

$$\sqrt{E\big[f^2(A_h)\big]}\sum\nolimits_{n=0}^{\infty}\frac{|\nu|^n}{n!}\cdot\sqrt{E\big[B_h^{2n}\big]}$$

being convergent for any complex number $\nu$. Implicit herein is that $f(A_h)$ is square integrable. This is implied by $A_h$ being square integrable. Any $n$–th moment in particular of $A_h$ has been computed in [**Y**, (4.d"), p.519] as:

$$n!\cdot\left(\frac{(-1)^n}{(n!)^2}+2\sum\nolimits_{k=0}^{n}\frac{(-1)^{n-k}}{(n-k)!\cdot(n+k)!}\cdot e^{-h\cdot\frac{k^2}{2}}\right).$$

In particular the second moment of $A_h$ is thus finite, as was to be shown. What regards the even order moments of the Brownian motion $B$ at time $h$, they are computed as:

$$E\big[B_h^{2n}\big]=\frac{(2h)^n}{n!}\cdot\Gamma\big(n+1/2\big),$$

for any non–negative integer $n$. The above series thus converges using the ratio test. This completes the proof of the normalized price being an entire function in the paramater $\nu$.

The proof of the six–term sum of §6 Theorem being a meromorphic function in $\nu$ is based on the Hermite functions being entire also in their degree. Indeed, from [**L**, (10.2.8), p.285] one has for any complex numbers $\mu$ and $z$ the representation:

$$H_\mu(z)=\frac{2^\mu\cdot\Gamma(1/2)}{\Gamma\big((1-\mu)/2\big)}\cdot\Phi\Big(-\frac{\mu}{2},\frac{1}{2},z^2\Big)+z\cdot\frac{2^\mu\cdot\Gamma(-1/2)}{\Gamma\big(-\mu/2\big)}\cdot\Phi\Big(\frac{1-\mu}{2},\frac{3}{2},z^2\Big).$$

Herein the reciprocal of the gamma function is entire by construction, and the confluent hypergeometric function $\Phi$ is entire in its first and third variable.

There is a *localization principle* for proving entire a function on the complex plane. Indeed, by definition one has to show analyticity at any fixed complex number. For this, it is sufficient to restrict the function to any relatively compact or compact neighborhood of this fixed complex number and prove the function thus obtained analytic.

What regards the hyperbolic terms, the claim is that as a function of $\nu$ they can be extended as analytic functions from non–negative $\nu$ to the whole complex plane. For this question it is sufficient to consider the function $A$ given by:

$$A(\nu)=\int_0^\infty H_{-(\nu+4)}\left(-\frac{\cosh(y+i\theta)}{\sqrt{2q}}\right)E_{b(\nu)}(h)(y+i\theta)\,dy\,,$$

where $b(\nu)$ equals $\pm\nu$ or $\pm(\nu+4)$ and recalling

$$E_\xi(h)(w) = e^{w\xi}\mathrm{Erfc}\left(\frac{w}{\sqrt{2h}} + \frac{\xi}{2}\sqrt{2h}\right),$$

for any complex number $w$. Extension of this function is essentially by reduction to the case where the real part of the degree $-(\nu+4)$ of the Hermite function factor in the integrand of $A(\nu)$ is negative, or equivalently, the real part of $\nu$ is bigger than minus four.

To fix ideas first consider the case where $\nu$ is such that the degree of the Hermite function factor in the integrand of $A(\nu)$ is a non–negative integer. This Hermite function is then the corresponding Hermite polynomial. In the integrand of $A(\nu)$ the absolute values of the Hermite function factor and the exponential function factor so have linear exponential order in the variable $y$. The decay to zero of the absolute value of the function $E_{b(\nu)}$, however, is of square exponential order in $y$. It thus dominates the asymptotic behaviour with $y$ to infinity. The absolute value of the integrand of $A(\nu)$ is so majorized by an integrable function, and $A$ can be extended to the above values of $\nu$.

For the general case of the reduction, fix any $\nu$ of non–positive real part that is not an integer. Apply the above localization principle and let $\nu$ belong to a sufficiently small compact neigborhood $U$ in the half–plane $\{\mathrm{Re}\,(z) \leq 0\}$. Shrinking $U$ if necessary assume that for any $\nu$ in $U$ the degree of the Hermite function in the integrand of $A(\nu)$ is not an integer. Using the recursion rule for Hermite functions of §2 express the Hermite function factor of $A(\nu)$ in terms of weighted Hermite functions of negative degrees. Further shrinking $U$ if necessary, assume that the so obtained relation represents $B$ on $U$. Herein the weighting factor for the respective Hermite functions are given as powers of $z = -(2q)^{-1/2}\cosh(y+i\theta)$ times polynomials in $\nu$. The absolute values of the polynomials in $\nu$ can be majoriozed uniformly on $U$. The problem thus reduces to majorize by an integrable function on the positive real line in the variable $y$ finitely many functions on $U$ times the positive real line sending $\nu$ and $y$ to:

$$H_{-(\nu+4+k)}\left(-\frac{\cosh(y+i\theta)}{\sqrt{2q}}\right) \cdot \cosh^\ell(y+i\theta) \cdot E_{b(\nu)}(h)(y+i\theta)\,,$$

where $k$, $\ell$ range over finitely many non–negative integers and $k$ is such that $\nu+4+k$ is positive. The leading terms of the asymptotic expansion of §2 for any Hermite function $H_\mu(z)$ with degree $\mu$ any complex number with negative real part has order $z^\mu$. Asymptotically with $y$ to infinity, the Hermite function with the smallest positive number $\nu+4+k$ thus dominates the other Hermite function factors in the above functions. With $\nu$ ranging over a compact set, there is a minimal such degree on $U$. Similarly, there are such majorizing choices $\ell^*$ for the factors $\cosh^\ell(y+i\theta)$ and $\nu^*$ for the absolute values of the factors $E_{b(\nu)}(h)(y+i\theta)$. A three–factor–majorizing function on the positive real line thus results whose asymptotic behaviour with $y$ to infinity is governed by the square–exponential decay to zero of the corresponding factor $E_{b(\nu^*)}(h)$ and which is integrable.

If the real part of $\nu$ is non–negative, the above argument holds in a simplified form. The upshot so is that any complex number not an integer has a sufficiently small compact neighborhood such that the absolute value of the integrand of $A(\nu)$ on $U$ times the positive real line can be majorized by an integrable function on the positive real line. Herein, compact neigborhoods can be replaced by relatively compact neigborhoods mutatis mutandis. Using Lebesgue Dominated Convergence, $A$ can thus be extended as a continuous function to the whole complex plane with the integers less than or equal to minus four deleted.

The idea for showing $A$ analytic as a function of $\nu$ on the complex plane with the integers less than or equal to minus four deleted is as follows. Show that differentiation of $A$ is by differentation under the integral sign and use that its integrand is entire as function of $\nu$. For this again first localize to $\nu$ in any sufficiently small compact neighborhood containing no integers less than or equal to minus four. The aim is then to majorize the absolute value of the derivative with respect to $\nu$ of the integrand of $A$ by an integrable function independent of $\nu$ as above. The above argument for getting such a majorizing function is based on a comparison of decay rates. The integrand of $A$ has one factor which on the positive real line decays to zero of square exponential order whereas the other factors explode of at most linear exponential order. This situation is preserved on differentiation with respect to the parameter $\nu$. In particular, differentiating with respect to the degree the asymptotic expansion for Hermite functions on the right half–plane gives an asymptotic expansion for this function's partial derivative with respect to the degree.

At this stage it remains to extend $A$ analytically to the intgers less than or equal to minus four. However, $A$ remains bounded in any punctured compact neigbourhood of such an integer. Thus $A$ can be extended to an entire function, completing the proof of §5 Theorem.

Proving the trigonometric terms entire is a direct application of the above localization principle to their following factors:

$$B(\nu) = \int_0^\theta H_{-(\nu+4)}\Big( - \frac{\cosh(\rho+i\phi)}{\sqrt{2q}}\Big) E_b(h)(\rho+i\phi)\Big) d\phi$$

where $b(\nu)$ equals $\nu$ or $\nu+2$. In fact, restrict $\nu$ to belong to any suffiently small compact neighborhood $U$ of any point fixed in the complex plane. The above integrand is analytic as a function in $\nu$ and smooth as a function in $\phi$. Thus the absolute values of its derivatives with respect to $\nu$ of any order are bounded on the product of $U$ and the closed interval between zero and $\theta$. In this local situation, using a standard consequence of Lebesgue Dominated Convergence, differentation of $B$ with respect to $\nu$ is thus by partial differentiation with respect to the parameter $\nu$ under the integral sign. This proves $B$ entire as a function in $\nu$.

**16. Trigonometric terms for unit radii:** At this point it remains to establish the exact form of the trigonometric terms $C_{\mathrm{trig},\theta,b(\nu)}(0)$ in §6 Remark for all $\nu$, thus extending §13 Lemma. This is essentially based on the above integral $B$ being entire. Indeed, $B$ being entire allows to first extend the complex identity underlying §6 Remark to all $\nu$. Thereupon taking real parts then completes the extension argument.

# Analytical ramifications of derivatives valuation: Asian options and special functions

by
**M. Schröder**
(Mannheim)

January 2002